\documentclass[12pt,leqno]{article}
\usepackage{amsmath, amssymb}

\newtheorem{theorem}{\bf Theorem}
\newtheorem{proposition}[theorem]{\bf Proposition}
\newtheorem{lemma}[theorem]{\bf Lemma}
\newtheorem{corollary}[theorem]{\bf Corollary}

\newtheorem{hypotheses}[theorem]{\bf Hypotheses}
\newtheorem{hypothesis}[theorem]{\bf Hypothesis}

\newcommand{\sect}[1]{\setcounter{equation}{0}\section{#1}}

\def\epsilon{\varepsilon}

\begin{document}

\Large \noindent 
{\bf Existence of traveling wave solutions}

\noindent
{\bf for a nonlocal bistable equation:} 

\noindent
{\bf an abstract approach}

\vspace*{0.8em}

\normalsize
\noindent Hiroki Yagisita

\noindent
Department of Mathematics, 
Faculty of Science, 
Kyoto Sangyo University

\noindent
Motoyama, Kamigamo, Kita-Ku, Kyoto-City, 603-8555, Japan

\vfill

\noindent {\bf Abstract} \ 
We consider traveling fronts to the nonlocal bistable equation  
\[u_t=\mu*u-u+f(u),\]
where $\mu$ is a Borel-measure on $\mathbb R$ with $\mu(\mathbb R)=1$ 
and $f$ satisfies $f(0)=f(1)=0$, $f<0$ in $(0,\alpha)$ 
and $f>0$ in $(\alpha,1)$ for some constant $\alpha \in (0,1)$. 
We do not assume that $\mu$ is absolutely continuous 
with respect to the Lebesgue measure. 
We show that there are a constant $c$ and a monotone function $\phi$ 
with $\phi(-\infty)=0$ and $\phi(+\infty)=1$ 
such that $u(t,x):=\phi(x+ct)$ is a solution 
to the equation, provided $f^\prime(\alpha)>0$. 
In order to prove this result, we would develop 
a recursive method for abstract monotone dynamical systems 
and apply it to the equation.

\vspace*{0.4em} 


\noindent
convolution model, integro-differential equation, 

\noindent
discrete bistable equation, nonlocal evolution equation.  

\vspace*{0.8em}


\vspace*{1.6em}

\noindent 
A proposed running title: Traveling Waves in Nonlocal Systems III.

\newpage




\sect{Introduction}
\( \, \, \, \, \, \, \, \) 
We would consider the following nonlocal analogue 
of bistable reaction-diffusion equations: 
\begin{equation}u_t=\mu*u-u+f(u).\end{equation} 
Here, $\mu$ is a Borel-measure on $\mathbb R$ with $\mu(\mathbb R)=1$ 
and the convolution is defined by 
\[(\mu*u)(x)=\int_{y\in\mathbb R}u(x-y)d\mu(y)\] 
for a bounded and Borel-measurable function $u$ on $\mathbb R$. 
The nonlinearity $f$ is a Lipschitz continuous function on $\mathbb R$ 
and satisfies $f(0)=f(\alpha)=f(1)=0$, $f<0$ in $(0,\alpha)$ 
and $f>0$ in $(\alpha,1)$ for some constant $\alpha \in (0,1)$. 
Then, $G(u):=\mu*u-u+f(u)$ is a map from the Banach space $L^\infty(\mathbb R)$ 
into $L^\infty(\mathbb R)$ and it is Lipschitz continuous. 
(We note that $u(x-y)$ is a Borel-measurable function on $\mathbb R^2$, 
and $\|u\|_{L^\infty(\mathbb R)}=0$ implies 
$\|\mu*u\|_{L^1(\mathbb R)}\leq 
\int_{y\in\mathbb R}(\int_{x\in\mathbb R}|u(x-y)|dx)d\mu(y)$=0.) 
So, because the standard theory of ordinary differential equations works, 
we have well-posedness of the equation (1.1) 
and it generates a flow in $L^\infty(\mathbb R)$. 

In this paper, we would show that there exists a traveling wave solution. 
The main result is the following: 
\begin{theorem} \ 
Suppose the bistable nonlinearity $f\in C^1(\mathbb R)$ satisfies 
\[f^\prime(\alpha)>0,\]
where $\alpha$ is the unique zero of $f$ in $(0,1)$. Then, 
there exist a constant $c$ and a monotone function $\phi$ on $\mathbb R$
with $\phi(-\infty)=0$ and $\phi(+\infty)=1$ 
such that $u(t,x):=\phi(x+ct)$ is a solution to (1.1). 
\end{theorem}
In this result, we do not assume that the measure $\mu$ is absolutely continuous 
with respect to the Lebesgue measure. For example, 
Theorem 1 can be applied to not only the integro-differential equation 
\[\frac{\partial u}{\partial t}(t,x)=\int_0^1u(t,x-y)dy-u(t,x)
-\lambda u(t,x)(u(t,x)-\alpha)(u(t,x)-1)\]
but also the discrete equation 
\[\frac{\partial u}{\partial t}(t,x)=u(t,x-1)-u(t,x)
-\lambda u(t,x)(u(t,x)-\alpha)(u(t,x)-1)\]
for all positive constants $\lambda$. 
In order to prove Theorem 1, we would develop a recursive method 
for abstract monotone dynamical systems and apply it to the semiflow 
generated by (1.1). It might be a generalization 
of the method of Remark 5.2 (4) in Chen [5].  

For the nonlocal bistable equation (1.1), 
Bates, Fife, Ren and Wang [4] obtained existence of traveling wave solutions, 
when the measure $\mu$ has a density function $J\in C^1(\mathbb R)$ 
with $J(y)=J(-y)$ and other little conditions for $\mu$ and $f$ hold. 
Chen [5] showed existence of traveling wave solutions, 
when it has a density function $J\in C^1(\mathbb R)$ 
and $f^\prime(u)<1$ and other little conditions hold. 
Recently, Coville [11] proved existence of traveling wave solutions, 
when it has a density function $J\in C(\mathbb R)$ 
and other little conditions hold. 

Bates, Fife, Ren and Wang [4] and 
Chen [5] studied uniqueness and stability of traveling wave solutions. 
Coville studied uniqueness and monotonicity of profiles of traveling waves 
in [10] and uniqueness of speeds [11]. 
Further, we note that the studies of [10, 11] are not limited 
when the nonlinearity is bistable but reach ignition, 
while our study is limited to bistable. 
In fact, his method of [11] is rather different from ours. 
See [9] on traveling wave solutions in bistable maps,
[2] time-periodic nonlocal bistable equations, 
[1] time-periodic bistable reaction-diffusion equations, 
e.g., [3, 6, 8, 14] discrete bistable equations, 
[7] nonlocal Burgers equations 
and [12, 13, 15] multistable reaction-diffusion equations.

In Section 2, we give abstract conditions and state that 
there exists a traveling wave solution provided the conditions. 
This result might generalize the method of Remark 5.2 (4) in Chen [5]. 
In Section 3, we prove abstract theorems mentioned in Section 2.
In Section 4, we show that the semiflow generated by (1.1) 
satisfies the conditions given in Section 2 
when $f^\prime(\alpha)>0$ and $\mu(\{0\})\not=1$ hold to prove Theorem 1. 
In Section 5, we recall some known results from [17]. The known results 
are used in Section 3.

\sect{Abstract theorems for monotone semiflows}
\( \, \, \, \, \, \, \, \) 
In this section, we would state some abstract results 
for existence of traveling waves in monotone semiflows. 
The results might generalize the method of Remark 5.2 (4) in Chen [5]. 
In the abstract, we would treat a bistable evolution system. Put a set 
of functions on $\mathbb R$; 
\[\mathcal M:=\{u\, |\, u \text{ is a monotone nondecreasing}\]
\[\text{ and left continuous function on } \mathbb R \text{ with } 0\leq u\leq 1\}.\] 

\vspace*{0.4em} 
%
%

The followings are our basic conditions for discrete dynamical systems: 
\begin{hypotheses} \ 
Let $Q_0$ be a map from $\mathcal M$ into $\mathcal M$. 

{\rm (i)} \ $Q_0$ is continuous in the following sense: 
If a sequence $\{u_k\}_{k\in \mathbb N}\subset \mathcal M$ 
converges to $u\in \mathcal M$ uniformly on every bounded interval, 
then the sequence $\{Q_0[u_k]\}_{k\in \mathbb N}$ converges to $Q_0[u]$ almost everywhere. 

{\rm (ii)} \ $Q_0$ is order preserving; i.e., 
\[u_1\leq u_2 \Longrightarrow Q_0[u_1]\leq Q_0[u_2]\] 
for all $u_1$ and $u_2\in \mathcal M$. 
Here, $u\leq v$ means that $u(x)\leq v(x)$ holds for all $x\in \mathbb R$. 

{\rm (iii)} \ $Q_0$ is translation invariant; i.e., 
\[T_{x_0}Q_0=Q_0T_{x_0}\] 
for all $x_0\in \mathbb R$. Here, $T_{x_0}$ is the translation operator 
defined by $(T_{x_0}[u])(\cdot):=u(\cdot-x_0)$. 

{\rm (iv)} \ $Q_0$ is bistable; i.e., there exists $\alpha \in (0,1)$ 
with $Q_0[\alpha]=\alpha$ such that
\[0<\gamma<\alpha \Longrightarrow Q_0[\gamma]<\gamma\] 
and 
\[\alpha<\gamma<1 \Longrightarrow \gamma<Q_0[\gamma]\] 
hold for all constant functions $\gamma$. 
\end{hypotheses}
{\bf Remark} \ If $Q_0$ satisfies Hypothesis 1 (iii), 
then $Q_0$ maps constant functions to constant functions. 

\vspace*{0.8em} 

\noindent
The following condition for discrete dynamical systems 
might be a little generalization of the condition 
in Remark 5.2 (4) of Chen [5]:  
\begin{hypothesis} \ 
Let $Q_0$ be a map from $\mathcal M$ into $\mathcal M$. 
If two constants $c_-$ and $c_+$ and two functions $\phi_-$ and $\phi_+\in\mathcal M$ 
satisfy $(Q_0[\phi_-])(x-c_-)\equiv\phi_-(x)$, $\phi_-(-\infty)=0$, $\phi_-(+\infty)=\alpha$, 
$(Q_0[\phi_+])(x-c_+)\equiv\phi_+(x)$, $\phi_+(-\infty)=\alpha$ and $\phi_+(+\infty)=1$, 
then the inequality $c_-<c_+$ holds.  
\end{hypothesis}

\vspace*{0.4em} 
 
The following states that existence of suitable {\it sub and super-solutions} 
implies existence of traveling wave solutions with an estimate of the speeds 
in the discrete dynamical systems on $\mathcal M$: 
\begin{theorem} \ 
Let a map $Q_0:\mathcal M \rightarrow \mathcal M$ satisfy {\rm Hypotheses 2} and {\rm 3}. 
Suppose a constant $\underline c$ and a function $\underline \psi\in\mathcal M$ 
with $\underline \psi(0)=0$ and $\underline \psi(+\infty)\in(\alpha,1]$ 
satisfy $\underline \psi(x)\leq (Q_0[\underline \psi])(x-\underline c)$. 
Suppose a constant $\overline c$ and a function $\overline \psi\in\mathcal M$ 
with $\overline \psi(-\infty)\in[0,\alpha)$ and $\overline \psi(0)=1$  
satisfy $(Q_0[\overline \psi])(x-\overline c)\leq \overline \psi(x)$. 
Then, there exist $c\in[\underline c,\overline c]$ 
and $\phi\in\mathcal M$ with $\phi(-\infty)=0$ 
and $\phi(+\infty)=1$ such that $(Q_0[\phi])(x-c)\equiv \phi(x)$ holds. 
\end{theorem}

\begin{corollary} \ 
Let a map $Q_0:\mathcal M \rightarrow \mathcal M$ satisfy {\rm Hypotheses 2} and {\rm 3}. 
Then, there exist $c\in\mathbb R$ and $\phi\in\mathcal M$ with $\phi(-\infty)=0$ 
and $\phi(+\infty)=1$ such that $(Q_0[\phi])(x-c)\equiv \phi(x)$ holds. 
\end{corollary}

\vspace*{0.4em} 

We add the following conditions 
to Hypotheses 2 for continuous dynamical systems on $\mathcal M$: 
\begin{hypotheses} \ 
Let $Q:=\{Q^t\}_{t\in[0,+\infty)}$ be a family of maps from $\mathcal M$ to $\mathcal M$. 

{\rm (i)} \ $Q$ is a semigroup; i.e., $Q^t\circ Q^s=Q^{t+s}$ 
for all $t$ and $s\in [0,+\infty)$. 

{\rm (ii)} \ $Q$ is continuous in the following sense: 
Suppose a sequence $\{t_k\}_{k\in \mathbb N}\subset[0,+\infty)$ converges to $0$, 
and $u\in \mathcal M$. Then, the sequence $\{Q^{t_k}[u]\}_{k\in \mathbb N}$ 
converges to $u$ almost everywhere. 
\end{hypotheses}

\noindent 
Instead of Hypothesis 3, we consider the following condition 
for continuous dynamical systems. It might also be a little generalization 
of the condition in Remark 5.2 (4) of Chen [5]:  
\begin{hypothesis} \ 
Let $Q:=\{Q^t\}_{t\in[0,+\infty)}$ be a family 
of maps from $\mathcal M$ to $\mathcal M$. 
If two constants $c_-$ and $c_+$ 
and two functions $\phi_-$ and $\phi_+\in\mathcal M$ 
with $\phi_-(-\infty)=0$, $\phi_-(+\infty)=\alpha$, 
$\phi_+(-\infty)=\alpha$ and $\phi_+(+\infty)=1$ 
satisfy $(Q^t[\phi_-])(x-c_-t)\equiv\phi_-(x)$
and $(Q^t[\phi_+])(x-c_+t)\equiv\phi_+(x)$ 
for all $t\in[0,+\infty)$, 
then the inequality $c_-<c_+$ holds.  
\end{hypothesis}
{\bf Remark} \ We could found similar hypotheses as Hypothesis 7 
for existence of traveling waves to reaction-diffusion equations 
with triple stable equilibria in [12, 13, 15].  

\vspace*{1.2em} 

As we would have Theorem 4 for the discrete dynamical systems, 
we would have the following for the continuous dynamical systems: 
\begin{theorem} \ 
Let $Q^t$ be a map from $\mathcal M$ to $\mathcal M$ for $t\in[0,+\infty)$. 
Suppose the map $Q^t$ satisfies {\rm Hypotheses 2} for all $t\in (0,+\infty)$, 
and the family $Q:=\{Q^t\}_{t\in[0,+\infty)}$ {\rm Hypotheses 6} and {\rm 7}. 
Then, the following holds {\rm :} 

Suppose a constant $\underline c$ 
and a function $\underline \psi\in\mathcal M$ 
with $\underline \psi(0)=0$ and $\underline \psi(+\infty)\in(\alpha,1]$ 
satisfy $\underline \psi(x)\leq (Q^t[\underline \psi])(x-\underline c t)$ 
for all $t\in[0,+\infty)$. 
Suppose a constant $\overline c$ and a function $\overline \psi\in\mathcal M$ 
with $\overline \psi(-\infty)\in[0,\alpha)$ and $\overline \psi(0)=1$ 
satisfy $(Q^t[\overline \psi])(x-\overline c t)\leq \overline \psi(x)$ 
for all $t\in[0,+\infty)$. 
Then, there exist $c\in[\underline c,\overline c]$ 
and $\phi\in\mathcal M$ with $\phi(-\infty)=0$ 
and $\phi(+\infty)=1$ such that $(Q^t[\phi])(x-ct)\equiv \phi(x)$ holds 
for all $t\in[0,+\infty)$. 
\end{theorem}

\sect{Proof of the abstract theorems}
\( \, \, \, \, \, \, \, \) 
In this section, we prove the theorems stated in Section 2 
by using known results recalled from [17] in Section 5. 

\vspace*{0.8em} 

\noindent
{\bf Proof of Theorem 4.}

[Step 0] \ In this step, we would give an intuitive explanation of our ideas. 
If you want to advance to exact proof at once, 
the step is recommended to be skipped.  

Because the map $Q_0 \, : \mathcal M \rightarrow \mathcal M$ is translation invariant, 
it is difficult to construct {\it traveling sub and super-solutions with the same speed} 
directly. So, we introduce a sequence 
of perturbed maps $Q_n \, : \mathcal M \rightarrow \mathcal M$ 
to break the translation invariance but to preserve the order. Then, 
we might construct sub and super-solutions $\underline \psi_n$ 
and $\overline \psi_n$ to the perturbed problem $Q_n[u]=u$ 
and also obtain a solution $\phi_n$ (i.e., $Q_n[\phi_n]=\phi_n$, 
$\phi_n(-\infty)=0$ and $\phi_n(+\infty)=1$) by order preserving property. 
In virtue of Hypothesis 3, we expect that the limit of a suitable subsequence 
of $(T_{-x_n}[\phi_n])(\cdot):=\phi_n(\cdot+x_n)$ 
solves the original problem. 

We shall explain more in detail but extremely inexactly. Let $n\in\mathbb N$.  
We put $\rho_n(x):=\frac{n+\frac{\overline c-\underline c}{2}}{n}
(x-\frac{\underline c+\overline c}{2})$. Then, 
the map $\, u \, \mapsto \, Q_n[u]:=Q_0[u\circ\rho_n]$ 
breaks the translation invariance but preserves the order. 
So, we may have a solution $\phi_n$ to $Q_n[\phi_n]=\phi_n$, $\phi_n(-\infty)=0$ 
and $\phi_n(+\infty)=1$. We take $y_n$ and $z_n$ such that $y_n\leq z_n$ 
and $0<\phi_n(y_n)<\alpha<\phi_n(z_n)<1$ hold. 

When a constant $c$ and a sequence $x_n$ satisfy $c=\frac{\underline c+\overline c}{2}
-\frac{\overline c-\underline c}{2}(\lim_{n\rightarrow\infty}\frac{x_n}{n})$, 
we see 
\[\lim_{n\rightarrow\infty}(\rho_n(x+x_n)-x_n)=x-c\] 
and, so, 
\[\lim_{n\rightarrow\infty}(T_{-x_n}\circ Q_n\circ T_{x_n})[u]
=\lim_{n\rightarrow\infty}T_{-x_n}[Q_0[(T_{x_n}[u])\circ\rho_n]]\] 
\[=\lim_{n\rightarrow\infty}Q_0[T_{-x_n}[(T_{x_n}[u])\circ\rho_n]
=Q_0[\lim_{n\rightarrow\infty}T_{-x_n}[(T_{x_n}[u])\circ\rho_n]]\] 
\[=Q_0[T_c[u]]=(Q_0\circ T_c)[u],\] 
where $(T_{x}[u])(\cdot):=u(\cdot-x)$. We might take a subsequence $n(k)$ 
such that there exist the limits $\phi_-:=\lim_{k\rightarrow\infty}T_{-y_n}[\phi_n]$, 
$\phi_+:=\lim_{k\rightarrow\infty}T_{-z_n}[\phi_n]$, 
$c_-:=\frac{\underline c+\overline c}{2}
-\frac{\overline c-\underline c}{2}(\lim_{k\rightarrow\infty}\frac{y_n}{n})$ 
and 
$c_+:=\frac{\underline c+\overline c}{2}
-\frac{\overline c-\underline c}{2}(\lim_{k\rightarrow\infty}\frac{z_n}{n})$. 

Therefore, we could expect that the two equalities 
\[(Q_0\circ T_{c_-})[\phi_-]
=\lim_{n\rightarrow\infty}(T_{-y_n}\circ Q_n\circ T_{y_n})[\phi_-]
=\lim_{k\rightarrow\infty}(T_{-y_n}\circ Q_n)[\phi_n]\]
\[=\lim_{k\rightarrow\infty}T_{-y_n}[Q_n[\phi_n]]
=\lim_{k\rightarrow\infty}T_{-y_n}[\phi_n]=\phi_-\]
and 
\[(Q_0\circ T_{c_+})[\phi_+]=\phi_+\] 
hold. In virtue of Hypothesis 3, the pair $(\phi_-,c_-)$ or $(\phi_+,c_+)$ 
might solve the original problem, as we obtain $c_+\leq c_-$ 
and $0<\phi_-(0)<\alpha<\phi_+(0)<1$.  

[Step 1] \ We show the inequality:  
\begin{equation}\underline c\leq\overline c.\end{equation} 

Suppose $\overline c<\underline c$. 
Then, there exists $N\in \mathbb N$ such that 
$\overline \psi (-\frac{\underline c-\overline c}{2}N)
<\alpha
<\underline \psi (+\frac{\underline c-\overline c}{2}N)$ 
holds. Hence, because 
$({Q_0}^N[\overline\psi])(x-\overline c N)\leq \overline\psi(x)$ 
and
$\underline\psi(x)\leq ({Q_0}^N[\underline\psi])(x-\underline c N)$ 
hold by Hypotheses 2 (ii) and (iii), we have 
$({Q_0}^N[\overline\psi])(-\frac{\underline c+\overline c}{2}N)
<\alpha 
<({Q_0}^N[\underline\psi])(-\frac{\underline c+\overline c}{2}N)$. 
It is a contradiction with Hypothesis 2 (ii). Therefore, 
(3.1) holds. 

[Step 2] \ We put a sequence $\{\rho_n\}_{n\in\mathbb N}$ 
of affine functions on $\mathbb R$ defined by 
\begin{equation}\rho_n(x):=\frac{n+\frac{\overline c-\underline c}{2}}{n}
(x-\frac{\underline c+\overline c}{2}).\end{equation} 
We define two sequences $\{A_n\}_{n\in\mathbb N}$ 
and $\{Q_n\}_{n\in\mathbb N}$ of maps from $\mathcal M$ 
to $\mathcal M$ by 
\[A_n[u]:=u\circ \rho_n\]
and 
\[Q_n:=Q_0\circ A_n.\] 
Then, the map $Q_n$ satisfies Hypothesis 2 (ii) 
for all $n\in\mathbb N$. 

[Step 3] \ 
We show the following: 
{\it Suppose a sequence $\{u_k\}_{k\in\mathbb N}\subset\mathcal M$ 
converges to $u\in\mathcal M$ almost everywhere. Then, 
$\lim_{k\rightarrow\infty}(Q_n[u_k])(x)=(Q_n[u])(x)$ 
holds for all $n\in\mathbb N$ and continuous points $x\in\mathbb R$
of $Q_n[u]$.} 

Let $n\in\mathbb N$. Then, the sequence $\{A_n[u_k]\}_{k\in\mathbb N}\subset\mathcal M$ 
converges to $A_n[u]\in\mathcal M$ almost everywhere. Hence, 
by Proposition 13, we have $\lim_{k\rightarrow\infty}(Q_0[A_n[u_k]])
$
$
(x)=(Q_0[A_n[u]])(x)$ 
for all continuous points $x\in\mathbb R$ of $Q_0[A_n[u]]$.

[Step 4] \ 
We take two sequences $\{\underline\psi_n\}_{n\in\mathbb N}$
and $\{\overline\psi_n\}_{n\in\mathbb N}\subset \mathcal M$ 
as 
\[\underline\psi_n(x):=\underline\psi(x-(n+\frac{\overline c-\underline c}{2}))\] 
and 
\[\overline\psi_n(x):=\overline\psi(x+(n+\frac{\overline c-\underline c}{2})).\] 
Then, we show 
$\underline\psi_n
\leq
{Q_n}^k[\underline\psi_n]
\leq
{Q_n}^{k+1}[\underline\psi_n] 
\leq
{Q_n}^{k+1}[\overline\psi_n]
\leq
{Q_n}^k[\overline\psi_n]
\leq
\overline\psi_n$ for all $k=0,1,2,\cdots$. 
Because $+n\leq x-\frac{\underline c+\overline c}{2}$ 
implies $x-\underline c\leq \rho_n(x)$ by (3.1), 
\begin{equation}\underline\psi_n(x-\underline c)
\leq (A_n[\underline\psi_n])(x)\end{equation} 
holds. Because $x-\frac{\underline c+\overline c}{2}\leq -n$ 
implies $\rho_n(x)\leq x-\overline c$ by (3.1), 
\begin{equation}(A_n[\overline\psi_n])(x)
\leq \overline\psi_n(x-\overline c)\end{equation} 
holds. From (3.3), (3.4) and $\underline\psi_n\leq\overline\psi_n$, we have 
$\underline\psi_n(x)\leq(Q_0[\underline\psi_n])(x-\underline c)
\leq (Q_n[\underline\psi_n])(x)
\leq (Q_n[\overline\psi_n])(x)
\leq (Q_0[\overline\psi_n])(x-\overline c)\leq\overline\psi_n(x)$. 
As $\underline\psi_n
\leq
{Q_n}^k[\underline\psi_n]
\leq
{Q_n}^{k+1}[\underline\psi_n] 
\leq
{Q_n}^{k+1}[\overline\psi_n]
\leq
{Q_n}^k[\overline\psi_n]
\leq
\overline\psi_n$ 
holds, 
$\underline\psi_n
\leq
Q_n[\underline\psi_n]
\leq
{Q_n}^{k+1}[\underline\psi_n]
\leq
{Q_n}^{k+2}[\underline\psi_n] 
\leq
{Q_n}^{k+2}[\overline\psi_n]
\leq
{Q_n}^{k+1}[\overline\psi_n]
\leq
Q_n[\overline\psi_n]
\leq
\overline\psi_n$ 
also holds. So, we have 
\[\underline\psi_n
\leq
{Q_n}^k[\underline\psi_n]
\leq
{Q_n}^{k+1}[\underline\psi_n] 
\leq
{Q_n}^{k+1}[\overline\psi_n]
\leq
{Q_n}^k[\overline\psi_n]
\leq
\overline\psi_n\] 
for all $n\in\mathbb N$ and $k=0,1,2,\cdots$. 
We put $\phi_n:=\lim_{k\rightarrow\infty}{Q_n}^k[\underline\psi_n]\in\mathcal M$. 
Then, 
\begin{equation}\underline\psi_n
\leq
\phi_n
\leq
\overline\psi_n\end{equation} 
holds for all $n\in\mathbb N$. By Step 3, we also have 
\begin{equation}Q_n[\phi_n]=\phi_n\end{equation} 
for all $n\in\mathbb N$. 

[Step 5] \ We take $N_0\in \mathbb N$ such that 
\begin{equation}0\leq \overline \psi (-N_0)
<\alpha
<\underline \psi (+N_0)\leq 1\end{equation} 
holds. Then, because 
$\phi_n(-(n+\frac{\overline c-\underline c}{2}+N_0))
\leq \overline\psi(-N_0)$ 
and 
$\underline\psi(+N_0) \leq 
\phi_n(+(n+\frac{\overline c-\underline c}{2}+N_0))$ 
hold from (3.5), for any $n\in\mathbb N$, there exist constants 
$y_n$ and $z_n$ such that 
\[\phi_n(y_n)\leq\frac{\overline\psi(-N_0)+\alpha}{2}
\leq\lim_{h\downarrow+0}\phi_n(y_n+h),\]
\[\phi_n(z_n)\leq\frac{\alpha+\underline\psi(+N_0)}{2}
\leq\lim_{h\downarrow+0}\phi_n(z_n+h)\]
and 
\begin{equation}-(n+\frac{\overline c-\underline c}{2}+N_0)
\leq y_n\leq z_n\leq 
+(n+\frac{\overline c-\underline c}{2}+N_0)\end{equation} 
hold. As we put functions 
\[\phi_{-,n}(\cdot):=\phi_n(\cdot+y_n)\in\mathcal M\] 
and 
\[\phi_{+,n}(\cdot):=\phi_n(\cdot+z_n)\in\mathcal M,\] 
we have 
\begin{equation}\phi_{-,n}(0)\leq\frac{\overline\psi(-N_0)+\alpha}{2}
\leq\lim_{h\downarrow+0}\phi_{-,n}(h)\end{equation} 
and 
\begin{equation}\phi_{+,n}(0)\leq\frac{\alpha+\underline\psi(+N_0)}{2}
\leq\lim_{h\downarrow+0}\phi_{+,n}(h).\end{equation} 
By Helly's theorem and (3.8), there exist a subsequence $\{n(k)\}_{k\in\mathbb N}
\subset\mathbb N$, two functions $\phi_-, \phi_+$, 
two constants $\xi_-$ and $\xi_+$ such that 
the two equalities 
\begin{equation}\phi_-(x)=\lim_{k\rightarrow\infty}\phi_{-,n(k)}(x)
\in\mathcal M\end{equation} 
and 
\[\phi_+(x)=\lim_{k\rightarrow\infty}\phi_{+,n(k)}(x)
\in\mathcal M\] hold 
almost everywhere in $x$ and the two equalities 
\begin{equation}\xi_-=\lim_{k\rightarrow\infty}\frac{y_{n(k)}}{n(k)}
\in[-1,+1]\end{equation} 
and   
\[\xi_+=\lim_{k\rightarrow\infty}\frac{z_{n(k)}}{n(k)}
\in[-1,+1]\] hold. 
From (3.9), (3.10) and (3.8), we have 
\begin{equation}\phi_-(0)\leq\frac{\overline\psi(-N_0)+\alpha}{2}
\leq\lim_{h\downarrow+0}\phi_-(h),\end{equation}
\begin{equation}\phi_+(0)\leq\frac{\alpha+\underline\psi(+N_0)}{2}
\leq\lim_{h\downarrow+0}\phi_+(h)\end{equation}
and 
\begin{equation}-1\leq\xi_-\leq \xi_+\leq+1.\end{equation} 

[Step 6] \ We show the following: 
{\it The two equalities 
\begin{equation}(Q_0[\phi_-])(x-c_-)\equiv\phi_-(x)\end{equation} 
and 
\begin{equation}(Q_0[\phi_+])(x-c_+)\equiv\phi_+(x)\end{equation} 
hold, where $c_-$ and $c_+$ are the constants defined by 
\[c_-:=\frac{\underline c+\overline c}{2}-\frac{\overline c-\underline c}{2}\xi_-\]
and
\[c_+:=\frac{\underline c+\overline c}{2}-\frac{\overline c-\underline c}{2}\xi_+.\]
Further, the inequality 
\begin{equation}\underline c\leq c_+\leq c_-\leq \overline c\end{equation}
holds. }

From (3.2) and (3.12), we see 
\[\lim_{k\rightarrow\infty}(\rho_{n(k)}(x+y_{n(k)})-y_{n(k)})=x-c_-\] 
for all $x\in\mathbb R$. 
Hence, by Lemma 14, (3.11) and 
$(A_n[\phi_n])(x+y_n)
=\phi_n(\rho_n(x+y_n))
=\phi_{-,n}(\rho_n(x+y_n)-y_n)$, we have 
\begin{equation}\lim_{k\rightarrow\infty}(A_{n(k)}[\phi_{n(k)}])(x+y_{n(k)})
=\phi_-(x-c_-)\end{equation} 
for all continuous points $x\in\mathbb R$ 
of $\phi_-(x-c_-)$. From (3.6), 
\begin{equation}\phi_{-,n}(x)=\phi_n(x+y_n)=(Q_n[\phi_n])(x+y_n)\end{equation}
\[=(Q_0[A_n[\phi_n]])(x+y_n)=(Q_0[(A_n[\phi_n])(\cdot+y_n)])(x)\] 
holds for all $n\in\mathbb N$ and $x\in\mathbb R$. 
By Proposition 13, (3.19), (3.20) and (3.11), we obtain 
\[(Q_0[\phi_-])(x-c_-)=(Q_0[\phi_-(\cdot-c_-)])(x)\]
\[=\lim_{k\rightarrow\infty}(Q_0[(A_{n(k)}[\phi_{n(k)}])(\cdot+y_{n(k)})])(x)\]
\[=\lim_{k\rightarrow\infty}\phi_{-,n(k)}(x)=\phi_-(x).\] 
Almost similarly as (3.16), we also obtain (3.17). 
Further, (3.18) follows from (3.1) and (3.15). 

[Step 7] \ By Proposition 13 and (3.16), we have 
\[Q_0[\phi_-(-\infty)]=(Q_0[\phi_-(-\infty)])(0)
=\lim_{k\rightarrow\infty}(Q_0[\phi_-(\cdot-k)])(0)\]
\[=\lim_{k\rightarrow\infty}(Q_0[\phi_-])(-k)
=(Q_0[\phi_-])(-\infty)=\phi_-(-\infty).\] 
Almost similarly, we also have 
$Q_0[\phi_-(+\infty)]=\phi_-(+\infty)$, 
$Q_0[\phi_+(-\infty)]=\phi_+(-\infty)$ 
and 
$Q_0[\phi_+(+\infty)]=\phi_+(+\infty)$ 
by Proposition 13, (3.16) and (3.17). From (3.7), (3.13) and (3.14), we see 
$0\leq \phi_-(-\infty)<\alpha$, $0<\phi_-(+\infty)\leq 1$, $0\leq\phi_+(-\infty)<1$ 
and $\alpha<\phi_+(+\infty)\leq 1$. Therefore, from Hypothesis 2 (iv), we obtain 
\begin{equation}\phi_-(-\infty)=0,\end{equation}
\begin{equation}\phi_-(+\infty)=\ \alpha \mbox{ or } 1,\end{equation}
\begin{equation}\phi_+(-\infty)=\ 0 \mbox{ or } \alpha\end{equation} 
and 
\begin{equation}\phi_+(+\infty)=1.\end{equation}

[Step 8] \ We show that $\phi_-(+\infty)\not=\alpha$ 
or $\phi_+(-\infty)\not=\alpha$ holds. Suppose 
that $\phi_-(+\infty)=\alpha$ 
and $\phi_+(-\infty)=\alpha$ hold. Then, from Hypothesis 3, (3.16), (3.17), 
(3.21) and (3.24), we have $c_-<c_+$. It is a contradiction 
with (3.18). So, we see that $\phi_-(+\infty)\not=\alpha$ 
or $\phi_+(-\infty)\not=\alpha$ holds. Hence, from (3.22) and (3.23), 
we also see that \[\phi_-(+\infty)=1 \ 
\mbox{ or } \ \phi_+(-\infty)=0\] holds. 
When $\phi_-(+\infty)=1$, we obtain the conclusion of Theorem 4 
with $c:=c_-$ and $\phi:=\phi_-$ because of (3.18), (3.21) and (3.16). 
When $\phi_+(-\infty)=0$, we obtain it 
with $c:=c_+$ and $\phi:=\phi_+$ because of (3.18), (3.24) and (3.17). 
\hfill 
$\blacksquare$ 

\vspace*{0.8em} 

\noindent
{\bf Proof of Corollary 5.}

We put functions $\underline \psi$ and $\overline \psi\in\mathcal M$ as 
\[\underline \psi(x)=0 \ \ \ (x\leq 0), 
\ \ \ \ \ \ \underline \psi(x)=\frac{\alpha+1}{2} \ \ \ (0<x)\]
and 
\[\overline \psi(x)=\frac{\alpha}{2} \ \ \ (x\leq -1), 
\ \ \ \ \ \ \overline \psi(x)=1 \ \ \ (-1<x).\]
Then, by Proposition 13 and Hypothesis 2 (iv), we have 
\[(Q_0[\underline \psi])(+\infty)
=\lim_{k\rightarrow\infty}(Q_0[\underline \psi])(k)
=\lim_{k\rightarrow\infty}(Q_0[\underline \psi(\cdot+k)])(0)\] 
\[=(Q_0[\underline \psi(+\infty)])(0)
=Q_0[\underline \psi(+\infty)]
>\underline \psi(+\infty).\] 
Almost similarly, we also have 
$(Q_0[\overline \psi])(-\infty)
<\overline \psi(-\infty)$. Hence, there exist constants 
$\underline c$ and $\overline c$ such that 
$\underline \psi(+\infty)\leq (Q_0[\underline \psi])(-\underline c)$
and $(Q_0[\overline \psi])(-1-\overline c)\leq \overline \psi(-\infty)$ hold. 
So, because $\underline \psi(x)\leq (Q_0[\underline \psi])(x-\underline c)$
and $(Q_0[\overline \psi])(x-\overline c)\leq \overline \psi(x)$ also hold 
for all $x\in\mathbb R$, in virtue of Theorem 4, we obtain the conclusion 
of Corollary 5.  
\hfill 
$\blacksquare$ 

%
%
%
%
\newpage 

\noindent
{\bf Proof of Theorem 8.}

[Step 1] \ 
By Lemma 17, 
the map $Q^t$ satisfies {\rm Hypothesis 3} for all $t\in (0,+\infty)$. 
So, by Theorem 4, for any $n\in\mathbb N$, 
there exist $c_n\in[\underline c,\overline c]$ and $\phi_n\in\mathcal M$ 
with $\phi_n(-\infty)=0$ and $\phi_n(+\infty)=1$ such that 
$(Q^{\frac{1}{2^n}}[\phi_n])(x-\frac{c_n}{2^n})
\equiv\phi_n(x)$ 
holds. Then, for any $n\in\mathbb N$, there exist constants 
$y_n$ and $z_n$ such that 
\[\phi_n(y_n)\leq\frac{\alpha}{2}
\leq\lim_{h\downarrow+0}\phi_n(y_n+h)\]
and 
\[\phi_n(z_n)\leq\frac{\alpha+1}{2}
\leq\lim_{h\downarrow+0}\phi_n(z_n+h)\]
hold. As we put functions 
\[\phi_{-,n}(\cdot):=\phi_n(\cdot+y_n)\in\mathcal M\] 
and 
\[\phi_{+,n}(\cdot):=\phi_n(\cdot+z_n)\in\mathcal M,\] 
we have 
\begin{equation}(Q^{\frac{1}{2^n}}[\phi_{-,n}])(x-\frac{c_n}{2^n})
\equiv\phi_{-,n}(x),\end{equation} 
\[(Q^{\frac{1}{2^n}}[\phi_{+,n}])(x-\frac{c_n}{2^n})
\equiv\phi_{+,n}(x),\] 
\begin{equation}\phi_{-,n}(0)\leq\frac{\alpha}{2}
\leq\lim_{h\downarrow+0}\phi_{-,n}(h)\end{equation} 
and 
\begin{equation}\phi_{+,n}(0)\leq\frac{\alpha+1}{2}
\leq\lim_{h\downarrow+0}\phi_{+,n}(h).\end{equation} 
By Helly's theorem, there exist a subsequence $\{n(k)\}_{k\in\mathbb N}
\subset\mathbb N$, two functions $\phi_-, \phi_+$ 
and a constant $c$ such that 
the two equalities 
\begin{equation}\phi_-(x)=\lim_{k\rightarrow\infty}\phi_{-,n(k)}(x)
\in\mathcal M\end{equation} 
and 
\[\phi_+(x)=\lim_{k\rightarrow\infty}\phi_{+,n(k)}(x)
\in\mathcal M\] 
hold almost everywhere in $x$ and the equality 
\begin{equation}c=\lim_{k\rightarrow\infty}c_{n(k)}
\in[\underline c,\overline c]\end{equation}
holds. From (3.26) and (3.27), we have 
\begin{equation}\phi_-(0)\leq\frac{\alpha}{2}
\leq\lim_{h\downarrow+0}\phi_-(h)\end{equation}
and 
\begin{equation}\phi_+(0)\leq\frac{\alpha+1}{2}
\leq\lim_{h\downarrow+0}\phi_+(h).\end{equation} 

[Step 2] \ We show the following: {\it The two equalities  
\begin{equation}(Q^t[\phi_-])(x-ct)\equiv\phi_-(x)\end{equation} 
and 
\begin{equation}(Q^t[\phi_+])(x-ct)\equiv\phi_+(x)\end{equation} 
hold for all $t\in[0,+\infty)$.}

Let $n_0\in\mathbb N$ and $m_0\in\mathbb N$. 
As $k\in\mathbb N$ is sufficiently large, 
\[(Q^{\frac{m_0}{2^{n_0}}}[\phi_{-,n(k)}])(x-c_{n(k)}\frac{m_0}{2^{n_0}})\] 
\[=((Q^{\frac{1}{2^{n(k)}}})^{m_02^{n(k)-n_0}}[\phi_{-,n(k)}])
(x-\frac{c_{n(k)}}{2^{n(k)}}m_02^{n(k)-n_0})
=\phi_{-,n(k)}(x)\] 
holds because of $n(k)\geq n_0$ and (3.25). Hence, by (3.28), (3.29), 
Lemma 14 and Proposition 13, we obtain 
\begin{equation}(Q^{\frac{m_0}{2^{n_0}}}[\phi_-])(x-c\frac{m_0}{2^{n_0}})
=\phi_-(x)\end{equation} 
for all $n_0\in\mathbb N$ and $m_0\in\mathbb N$. 

Let $t\in[0,+\infty)$. Then, by (3.34), 
there exists a sequence $\{t_k\}_{k\in\mathbb N}\subset[0,+\infty)$ 
with $\lim_{k\rightarrow\infty}t_k=0$ such that 
$(Q^{t+t_k}[\phi_-])(x-c(t+t_k))=\phi_-(x)$ 
holds for all $k\in\mathbb N$. So, by 
$(Q^{t_k}[(Q^t[\phi_-])(\cdot-ct)])(x-ct_k)=(Q^{t+t_k}[\phi_-])(x-c(t+t_k))$ 
and Lemma 15, we obtain $(Q^t[\phi_-])(x-ct)=\phi_-(x)$.  

Almost similarly as (3.32), we also obtain (3.33). 

[Step 3] \ By Proposition 13 and (3.32), we have 
\[Q^t[\phi_-(-\infty)]=(Q^t[\phi_-(-\infty)])(0)
=\lim_{k\rightarrow\infty}(Q^t[\phi_-(\cdot-k)])(0)\]
\[=\lim_{k\rightarrow\infty}(Q^t[\phi_-])(-k)
=(Q^t[\phi_-])(-\infty)=\phi_-(-\infty).\] 
Almost similarly, we also have 
$Q^t[\phi_-(+\infty)]=\phi_-(+\infty)$, 
$Q^t[\phi_+(-\infty)]=\phi_+(-\infty)$ 
and 
$Q^t[\phi_+(+\infty)]=\phi_+(+\infty)$ 
by Proposition 13, (3.32) and (3.33). From (3.30) and (3.31), we see 
$0\leq \phi_-(-\infty)<\alpha$, $0<\phi_-(+\infty)\leq 1$, $0\leq\phi_+(-\infty)<1$ 
and $\alpha<\phi_+(+\infty)\leq 1$. Therefore, from Hypothesis 2 (iv), we obtain 
\begin{equation}\phi_-(-\infty)=0,\end{equation}
\begin{equation}\phi_-(+\infty)=\ \alpha \mbox{ or } 1,\end{equation}
\begin{equation}\phi_+(-\infty)=\ 0 \mbox{ or } \alpha\end{equation} 
and 
\begin{equation}\phi_+(+\infty)=1.\end{equation}

[Step 4] \ We show that $\phi_-(+\infty)\not=\alpha$ 
or $\phi_+(-\infty)\not=\alpha$ holds. Suppose 
that $\phi_-(+\infty)=\alpha$ 
and $\phi_+(-\infty)=\alpha$ hold. Then, from Hypothesis 3, (3.32), (3.33), 
(3.35) and (3.38), we have the contradiction $c<c$. So, we see that $\phi_-(+\infty)\not=\alpha$ or $\phi_+(-\infty)\not=\alpha$ holds. 
Hence, from (3.36) and (3.37), 
we also see that \[\phi_-(+\infty)=1 \ 
\mbox{ or } \ \phi_+(-\infty)=0\] holds. 
When $\phi_-(+\infty)=1$, we obtain the conclusion of Theorem 8 
with $\phi:=\phi_-$. When $\phi_+(-\infty)=0$, we obtain it 
with $\phi:=\phi_+$. 
\hfill 
$\blacksquare$ 

\sect{Proof of Theorem 1} 
\( \, \, \, \, \, \, \, \) 
We recall that $\mu$ is a Borel-measure on $\mathbb R$ 
with $\mu(\mathbb R)=1$, $f$ is a Lipschitz continuous function on $\mathbb R$ 
and satisfies $f(0)=f(\alpha)=f(1)=0$, $f<0$ in $(0,\alpha)$ 
and $f>0$ in $(\alpha,1)$ for some constant $\alpha \in (0,1)$  
and the set $\mathcal M$ has been defined at the beginning of Section 2. 
Then, in virtue of Lemma 7 of [18], Lemma 8 of [18] and Proposition 10 of [18], 
$Q^t \ (t\in(0,+\infty))$ satisfies Hypotheses 2 and $Q$ Hypotheses 6 
for the semiflow $Q=\{Q^t\}_{t\in[0,+\infty)}$ on $\mathcal M$ 
generated by (1.1). So, if we would confirm that this semiflow on $\mathcal M$ 
satisfies Hypothesis 7, then we could make Theorem 8 of Section 2 work. 
In this section, we confirm it when $f^\prime(\alpha)>0$ 
and $\mu(\{0\})\not=1$ hold and construct sub and super-solutions 
to prove Theorem 1.

\vspace*{0.2em} 

First, we consider the linear equation 
\begin{equation}v_t=\hat \mu*v.\end{equation} 
It generates a flow on the Banach space $BC(\mathbb R)$
when $\hat\mu(\mathbb R)<+\infty$. 
Here, $BC(\mathbb R)$ denote the set of bounded and continuous functions on $\mathbb R$. 
We have the following for this flow on $BC(\mathbb R)$: 
\begin{proposition} \ 
Let $\hat \mu$ be a Borel-measure on $\mathbb R$ with $\hat\mu(\mathbb R)<+\infty$. 
Let $\, \hat P:\, BC(\mathbb R)\rightarrow BC(\mathbb R)\, $ be the time $1$ map 
of the flow on $BC(\mathbb R)$ generated by the linear equation (4.1). 
Then, there exists a Borel-measure $\hat \nu$ on $\mathbb R$ 
with $\hat \nu(\mathbb R)<+\infty$ such that 
\[\hat P[v]=\hat \nu*v\] 
holds for all $v\in BC(\mathbb R)$. Further, the equality  
\begin{equation}\log \int_{y\in\mathbb R}e^{\lambda y}d\hat \nu(y) 
=\int_{y\in \mathbb R}e^{\lambda y}d\hat\mu(y)\end{equation}
holds for all $\lambda \in\mathbb R$. 
\end{proposition}
{\bf Proof.} \ 
From Lemma 16 of [18], 
there exists a Borel-measure $\hat \nu$ on $\mathbb R$ 
with $\hat \nu(\mathbb R)<+\infty$ such that 
\begin{equation}\hat P[v]=\hat \nu*v\end{equation} 
holds for all $v\in BC(\mathbb R)$. Further, from Lemma 16 of [18], 
if $v$ is a nonnegative, bounded and continuous function on $\mathbb R$, 
then the inequality  
\[\hat\mu*v\leq\hat\nu*v\] 
holds. So, because 
\[\int_{y\in \mathbb R}e^{\lambda y}d\hat\mu(y)
=\lim_{n\rightarrow\infty}\int_{y\in \mathbb R}\min\{e^{\lambda y},n\}d\hat\mu(y)
=\lim_{n\rightarrow\infty}(\hat\mu*\min\{e^{-\lambda x},n\})(0)
\]
\[\leq\lim_{n\rightarrow\infty}(\hat\nu*\min\{e^{-\lambda x},n\})(0)
=\lim_{n\rightarrow\infty}\int_{y\in \mathbb R}\min\{e^{\lambda y},n\}d\hat\nu(y)
=\int_{y\in\mathbb R}e^{\lambda y}d\hat\nu(y)\] 
holds, $\int_{y\in \mathbb R}e^{\lambda y}d\hat\mu(y)=+\infty$
implies 
$\int_{y\in\mathbb R}e^{\lambda y}d\hat\nu(y)=+\infty$.  
Therefore, it is sufficient if we show that 
the equality (4.2) holds when 
\begin{equation}\int_{y\in \mathbb R}
e^{\lambda y}d\hat\mu(y)<+\infty.\end{equation}
Let $\lambda\in\mathbb R$. Suppose (4.4). 

Let $X_\lambda$ denote 
the set of continuous functions $u$ on $\mathbb R$ 
with $\sup_{x\in \mathbb R}\frac{|u(x)|}{1+e^{-\lambda x}}<+\infty$. 
Then, $X_\lambda$ is a Banach space with the norm 
$\|u\|_{X_\lambda}:=\sup_{x\in \mathbb R}\frac{|u(x)|}{1+e^{-\lambda x}}$. 
Let $u\in X_\lambda$. Then, for any $x$ and $y\in\mathbb R$, we have 
\[\sup_{h\in[-1,+1]}|u((x+h)-y)-u(x-y)|\]
\[\leq
\|u\|_{X_\lambda}\sup_{h\in[-1,+1]}((1+e^{-\lambda((x+h)-y)})+(1+e^{-\lambda(x-y)}))
\] 
\[\leq
\|u\|_{X_\lambda}(\sup_{h\in[-1,+1]}
((1+e^{-\lambda(x+h)})+(1+e^{-\lambda x})))(1+e^{\lambda y}).
\]  
Hence, from (4.4), the function $\hat\mu*u$ is continuous. Because 
\[\sup_{x\in \mathbb R}\frac{|(\hat\mu*u)(x)|}{1+e^{-\lambda x}}
\leq \sup_{x\in \mathbb R}\int_{y\in\mathbb R}
\frac{|u(x-y)|}{1+e^{-\lambda (x-y)}}(1+e^{\lambda y})d\hat\mu(y)\]
\[\leq \left(\int_{y\in\mathbb R}(1+e^{\lambda y})d\hat\mu(y)\right) 
\|u\|_{X_\lambda}\]
also holds, the map $ \, u \, \mapsto \, \hat\mu*u \, $ 
is a bounded and linear operator in the Banach space $X_\lambda$.  
Let $\, \hat P_\lambda: \, X_\lambda\rightarrow X_\lambda \,$ be the time $1$ map 
of the flow on $X_\lambda$ generated by the linear equation (4.1). 

Suppose $\lambda>0$. Let $\bar \lambda\in (0,\lambda)$. Then, we see 
\begin{equation}\lim_{n\rightarrow\infty}
\|\min\{e^{-\bar \lambda x},n\}-e^{-\bar \lambda x}\|_{X_\lambda}\end{equation}
\[\leq \lim_{n\rightarrow\infty}\sup_{x\in (-\infty,-\frac{1}{\bar \lambda}\log n)}
\frac{e^{-\bar \lambda x}}{1+e^{-\lambda x}}\]
\[\leq \lim_{n\rightarrow\infty}\sup_{x\in (-\infty,-\frac{1}{\bar \lambda}\log n)}
e^{(\lambda-\bar \lambda) x}=0.\]
The function $v(t,x)
:=e^{(\int_{y\in \mathbb R}e^{\bar \lambda y}d\hat\mu(y))t-\bar \lambda x}$ 
is a solution to (4.1) in the phase space $X_\lambda$. 
Hence, by (4.3) and (4.5), 
\[\int_{y\in\mathbb R}e^{\bar \lambda y}d\hat\nu(y)
=\lim_{n\rightarrow\infty}\int_{y\in\mathbb R}\min\{e^{\bar \lambda y},n\}d\hat\nu(y)\] 
\[=\lim_{n\rightarrow\infty}(\hat\nu*\min\{e^{-\bar \lambda x},n\})(0)
=\lim_{n\rightarrow\infty}(\hat P[\min\{e^{-\bar \lambda x},n\}])(0)\] 
\[=\lim_{n\rightarrow\infty}(\hat P_\lambda[\min\{e^{-\bar \lambda x},n\}])(0)
=(\hat P_\lambda[e^{-\bar \lambda x}])(0)
=e^{\int_{y\in \mathbb R}e^{\bar \lambda y}d\hat\mu(y)}\] 
holds for all $\bar \lambda\in (0,\lambda)$. So, we have 
\[\int_{y\in\mathbb R}e^{\lambda y}d\hat\nu(y)
=\lim_{\bar \lambda \, \uparrow \, \lambda}\int_{y\in\mathbb R}e^{\bar \lambda y}d\hat\nu(y)
=\lim_{\bar \lambda \, \uparrow \, \lambda}
e^{\int_{y\in \mathbb R}e^{\bar \lambda y}d\hat\mu(y)}
=e^{\int_{y\in \mathbb R}e^{\lambda y}d\hat\mu(y)}.\]

When $\lambda<0$, we could also prove it almost similarly as $\lambda>0$.  

Because 
$e^{(\int_{y\in \mathbb R}1d\hat\mu(y))t}$ 
is a solution to (4.1), from (4.3), we see 
\[\int_{y\in\mathbb R}1d\hat\nu(y)=(\hat\nu*1)(0)=(\hat P[1])(0)
=e^{\int_{y\in \mathbb R}1d\hat\mu(y)}.\]
So, the equality (4.2) also holds when $\lambda=0$. 
\hfill 
$\blacksquare$ 

\vspace*{0.4em} 

In [18], the author has recalled a method to estimate the spreading speeds 
in monostable systems by Weinberger [16]. 
Combining Proposition 9 with the method, 
we have the following: 
\begin{lemma} \ 
Suppose a constant $\sigma$ satisfies $0<\sigma<f^\prime(\alpha)$. 
Then, the following two hold {\rm :} 

{\rm (i)} \ Let $c_-\in\mathbb R$. Let $\phi_-$ be a monotone function 
on $\mathbb R$ with $\phi_-(-\infty)=0$ and $\phi_-(+\infty)=\alpha$. 
Suppose $u_-(t,x):=\phi_-(x+c_-t)$ is a solution to (1.1). Then, 
\[\inf_{\lambda_->0}
\frac{\int_{y\in\mathbb R}e^{\lambda_- y}d\mu(y)-1+\sigma}{\lambda_-}
\leq -c_-\]
holds. 

{\rm (ii)} \ Let $c_+\in\mathbb R$. Let $\phi_+$ be a monotone function 
on $\mathbb R$ with $\phi_+(-\infty)=\alpha$ and $\phi_+(+\infty)=1$. 
Suppose $u_+(t,x):=\phi_+(x+c_+t)$ is a solution to (1.1). Then, 
\[\inf_{\lambda_+>0}
\frac{\int_{y\in\mathbb R}e^{-\lambda_+ y}d\mu(y)-1+\sigma}{\lambda_+}
\leq c_+\]
holds.  
\end{lemma}
{\bf Proof.} \ 
[Step 1] \ In this step, we show (i). 

We put a Borel-measure $\hat\mu:=\mu$ 
and a Lipschitz continuous function $\hat f(u):=-\frac{1}{\alpha}f(-\alpha(u-1))$. 
Then, we see $\hat f(0)=\hat f(1)=0$, $f>0$ in $(0,1)$ and 
\begin{equation}0<\sigma<f^\prime(\alpha)={\hat f}^\prime(0).\end{equation} 
Further, we put a monotone function $\hat \phi(z):=-\frac{1}{\alpha}\phi_-(z)+1$ 
with $\hat \phi(-\infty)=1$ and $\hat \phi(+\infty)=0$. 
Then, the function $u(t,x):=\hat\phi(x+c_-t)$ is a solution to 
\begin{equation}u_t=\hat\mu*u-u+\hat f(u).\end{equation} 

Let $\, \hat P:\, BC(\mathbb R)\rightarrow BC(\mathbb R)\, $ be the time $1$ map 
of the flow on $BC(\mathbb R)$ generated by the linear equation (4.1). 
Then, by Proposition 9, there exists a Borel-measure $\hat \nu$ on $\mathbb R$ 
with $\hat \nu(\mathbb R)<+\infty$ such that 
\begin{equation}\hat P[v]=\hat \nu*v\end{equation} 
holds for all $v\in BC(\mathbb R)$. Further, the equality  
\begin{equation}\log \int_{y\in\mathbb R}e^{\lambda y}d\hat \nu(y) 
=\int_{y\in \mathbb R}e^{\lambda y}d\hat\mu(y)\end{equation}
holds for all $\lambda \in\mathbb R$. 
Let $\, \tilde P:\, BC(\mathbb R)\rightarrow BC(\mathbb R)\, $ be the time $1$ map 
of the flow on $BC(\mathbb R)$ generated by the linear equation 
\[v_t=\hat\mu*v-v+\sigma v.\] 
Then, from (4.8) and (4.9), as $\tilde \nu$ is the Borel-measure on $\mathbb R$ 
defined by 
\[\tilde \nu:=e^{-1+\sigma}\hat \nu,\] 
\begin{equation}\tilde P[v]=\tilde \nu*v\end{equation} 
holds for all $v\in BC(\mathbb R)$ and   
\begin{equation}\log \int_{y\in\mathbb R}e^{\lambda y}d\tilde \nu(y) 
=\int_{y\in \mathbb R}e^{\lambda y}d\hat\mu(y)-1+\sigma\end{equation}
holds for all $\lambda \in\mathbb R$. 
Because $\tilde\nu(\mathbb R)=(\tilde\nu*1)(0)=(\tilde P[1])(0)=e^\sigma$ 
holds from (4.10) and $\hat\mu(\mathbb R)=1$, we also have 
\begin{equation}1<\tilde \nu(\mathbb R)<+\infty.\end{equation} 

Let $\, \tilde Q_0:\, \mathcal B\rightarrow \mathcal B\, $ be the time $1$ map 
of the semiflow on $\mathcal B$ generated by the equation (4.7), 
where $\mathcal B$ is the set defined by 
\[\mathcal B:=\{u \, | \, 
u \mbox{ is a continuous function on } \mathbb R \mbox{ with } 0\leq u\leq 1\}.\] 
Then, from Lemma 17 of [18], (4.6) and (4.10), 
there exists a positive constant $\varepsilon$ such that the inequality   
\[\tilde \nu*u=\tilde P[u]\leq \tilde Q_0[u]\] 
holds for all $u\in\mathcal B$ with $u\leq \varepsilon$. 
Therefore, by Proposition 15 of [18] and (4.12), we obtain the inequality 
\[\inf_{\lambda>0}\frac{1}{\lambda}\log\int_{y\in\mathbb R}
e^{\lambda y}d\tilde \nu(y)
\,
\leq 
\,
-c_-.\] 
So, from (4.11), we obtain 
\begin{equation}\inf_{\lambda>0}
\frac{\int_{y\in \mathbb R}e^{\lambda y}d\hat\mu(y)-1+\sigma}{\lambda}
\,
\leq 
\,
-c_-.\end{equation}

[Step 2] \ We show (ii). 
Let $\hat\mu$ be the Borel-measure on $\mathbb R$ such that 
\begin{equation}\hat\mu((-\infty,y))=\mu((-y,+\infty))\end{equation} 
holds for all $y\in\mathbb R$. 
We put a Lipschitz continuous function 
$\hat f(u):=\frac{1}{1-\alpha}f((1-\alpha)u+\alpha)$ 
and a monotone function 
$\hat \phi(z):=\frac{1}{1-\alpha}(\phi_+(-z)-\alpha)$ 
with $\hat \phi(-\infty)=1$ and $\hat \phi(+\infty)=0$. 
Then, the function $u(t,x):=\hat\phi(x-c_+t)$ is a solution to 
\[u_t=\hat\mu*u-u+\hat f(u).\]  
So, almost similarly as (4.13), we also obtain 
\[\inf_{\lambda>0}
\frac{\int_{y\in \mathbb R}e^{\lambda y}d\hat\mu(y)-1+\sigma}{\lambda}
\,
\leq 
\,
c_+.\] 
Hence, the conclusion of (ii) follows from (4.14). 
\hfill 
$\blacksquare$ 

\vspace*{0.4em} 

\begin{lemma} \ 
Suppose $\mu(\{0\})\not=1$. Then, 
\begin{equation}0<\inf_{\lambda_->0}
\frac{\int_{y\in\mathbb R}(e^{\lambda_- y}-1)d\mu(y)+\sigma}{\lambda_-}
+\inf_{\lambda_+>0}
\frac{\int_{y\in\mathbb R}(e^{-\lambda_+ y}-1)d\mu(y)+\sigma}{\lambda_+}
\end{equation} 
holds for all $\sigma\in(0,+\infty)$. 
\end{lemma}
{\bf Proof.} \ Because 
$\inf_{\lambda_->0}\int_{y\in\mathbb R}e^{\lambda_- y}d\mu(y)\not=+\infty$ 
implies 
$\lim_{\lambda_-\downarrow+0}\int_{y\in\mathbb R}e^{\lambda_- y}
$
$
d\mu(y)=1$, 
we see $\lim_{\lambda_-\downarrow+0}
\frac{\int_{y\in\mathbb R}(e^{\lambda_- y}-1)d\mu(y)+\sigma}{\lambda_-}
=+\infty$. Similarly, we also see 
$\lim_{\lambda_+\downarrow+0}
\frac{\int_{y\in\mathbb R}(e^{-\lambda_+ y}-1)d\mu(y)+\sigma}{\lambda_+}
=+\infty$. So, $\inf_{\lambda_->0}
\frac{\int_{y\in\mathbb R}(e^{\lambda_- y}-1)d\mu(y)+\sigma}{\lambda_-}
\not=-\infty$ and $\inf_{\lambda_+>0}
\frac{\int_{y\in\mathbb R}(e^{-\lambda_+ y}-1)d\mu(y)+\sigma}{\lambda_+}
\not=-\infty$ hold. 
Hence, it is sufficient if we show that 
the inequality (4.15) holds 
when \begin{equation}\inf_{\lambda_->0}
\frac{\int_{y\in\mathbb R}(e^{\lambda_- y}-1)d\mu(y)+\sigma}{\lambda_-}
\not=+\infty\end{equation} and \begin{equation}\inf_{\lambda_+>0}
\frac{\int_{y\in\mathbb R}(e^{-\lambda_+ y}-1)d\mu(y)+\sigma}{\lambda_+}
\not=+\infty\end{equation} hold. 
Suppose (4.16) and (4.17) hold. 
Then, we have 
\begin{equation}\int_{y\in\mathbb R}|y|d\mu(y)<+\infty.\end{equation} 
Because of $\mu(\{0\})\not=\mu(\mathbb R)$, we see that  
$\mu((-\infty,0))\not=0$ or $\mu((0,+\infty))\not=0$ holds. 

Suppose $\mu((-\infty,0))\not=0$. Then, 
$\lim_{\lambda_+\rightarrow+\infty}
\frac{\int_{y\in\mathbb R}(e^{-\lambda_+ y}-1)d\mu(y)+\sigma}{\lambda_+}
=+\infty$ holds. Hence, there exists $N\in\mathbb N$ such that 
$\inf_{\lambda_+>0}
\frac{\int_{y\in\mathbb R}(e^{-\lambda_+ y}-1)d\mu(y)+\sigma}{\lambda_+}
=
%
%
%
%
%
%
\inf_{\lambda_+\in(0,N]}
\frac{\int_{y\in\mathbb R}(e^{-\lambda_+ y}-1)d\mu(y)+\sigma}{\lambda_+}$ 
holds. So, from (4.18), we have 
\[\inf_{\lambda_->0}
\frac{\int_{y\in\mathbb R}(e^{\lambda_- y}-1)d\mu(y)+\sigma}{\lambda_-}
+\inf_{\lambda_+>0}
\frac{\int_{y\in\mathbb R}(e^{-\lambda_+ y}-1)d\mu(y)+\sigma}{\lambda_+}\] 
\[\geq \inf_{\lambda_->0}
(\int_{y\in\mathbb R}yd\mu(y)+\frac{\sigma}{\lambda_-})
+\inf_{\lambda_+\in(0,N]}
(-\int_{y\in\mathbb R}yd\mu(y)+\frac{\sigma}{\lambda_+})
=\frac{\sigma}{N}.\] 

When $\mu((0,+\infty))\not=0$, we also have the inequality (4.15) 
almost similarly as $\mu((-\infty,0))\not=0$. 
\hfill 
$\blacksquare$

\vspace*{0.8em}

The following gives sub and super-solutions: 
\begin{lemma} \ 
Let a Lipschitz continuous function $\hat f$ on $\mathbb R$ satisfy  
\begin{equation}\hat f>0 \ \ \ (u<0), \ \ \ 
\hat f=f \ \ \ (0\leq u\leq 1), \ \ \ 
\hat f<0 \ \ \ (1<u).\end{equation}
Let a function $\rho\in C^1(\mathbb R)$ satisfy 
$\rho=0$ in $(-\infty,0]$, $\rho=1$ in $[1,+\infty)$ 
and $\rho^\prime>0$ in $(0,1)$. 
Suppose a positive constant $\varepsilon$ is sufficiently small. 
Then, the function $\underline u(t,x)
:=\rho(\varepsilon x-\frac{t}{\varepsilon})-\frac{1-\alpha}{4}$ 
is a sub-solution to the equation 
\begin{equation}u_t=\mu*u-u+\hat f(u)\end{equation} 
and  
the function $\overline u(t,x)
:=\rho(\varepsilon x+\frac{t}{\varepsilon}+1)+\frac{\alpha}{4}$ 
is a super-solution to (4.20). 
\end{lemma}
{\bf Proof.} \ We put a positive constant $\delta$ as 
\[\delta:=\min\left\{
\min_{u\in[-\frac{1-\alpha}{4},-\frac{1-\alpha}{8}]
\cup[1-\frac{1-\alpha}{2},1-\frac{1-\alpha}{4}]}\hat f(u), 
\, 
\min_{u\in[\frac{\alpha}{4},\frac{\alpha}{2}]
\cup[1+\frac{\alpha}{8},1+\frac{\alpha}{4}]}(-\hat f(u))
\right\}.\]
We also put a positive constant $C$ as 
\[C:=\min\left\{\rho^\prime(z) \, | \, 
\min\{\frac{1-\alpha}{8},\frac{\alpha}{4}\}
\leq \rho(z) \leq 
1-\min\{\frac{1-\alpha}{4},\frac{\alpha}{8}\}
\right\}.\] 
Then, we see 
\[\frac{1}{\varepsilon}\rho^\prime(z)+\hat f\left(\rho(z)-\frac{1-\alpha}{4}\right)
\geq \min\left\{\frac{1}{\varepsilon}C+
\min_{u\in[-\frac{1-\alpha}{8},1-\frac{1-\alpha}{2}]}
\hat f(u), \, \delta\right\}=\delta\]
for all $z\in\mathbb R$. 
We also see 
\[\frac{1}{\varepsilon}\rho^\prime(z)-\hat f\left(\rho(z)+\frac{\alpha}{4}\right)
\geq \min\left\{\frac{1}{\varepsilon}C+
\min_{u\in[\frac{\alpha}{2},1+\frac{\alpha}{8}]}
(-\hat f(u)), \, \delta\right\}=\delta\]
for all $z\in\mathbb R$. There exists $N\in\mathbb N$ such that 
$\mu(\mathbb R\setminus (-N,+N))<\delta$ holds. So, 
\[\left|\int_{y\in\mathbb R}\rho(\varepsilon (x-y)+s)d\mu(y)-\rho(\varepsilon x+s)\right|\]
\[\leq 
\int_{y\in\mathbb R}|\rho(\varepsilon (x-y)+s)-\rho(\varepsilon x+s)|d\mu(y)\]
\[\leq \int_{y\in(-N,+N)}|\rho(\varepsilon (x-y)+s)-\rho(\varepsilon x+s)|d\mu(y)
+\mu(\mathbb R\setminus (-N,+N))\]
\[\leq \sup_{h\in(-\varepsilon N,+\varepsilon N), \, z\in\mathbb R}|\rho(z+h)-\rho(z)|
+\mu(\mathbb R\setminus (-N,+N))\leq \delta\] 
holds for all $s$ and $x\in\mathbb R$. Therefore, we have 
\[-(\mu*\underline u-\underline u)\leq \delta \leq -{\underline u}_t+\hat f(\underline u)\] 
and 
\[\mu*\overline u-\overline u\leq \delta \leq {\overline u}_t-\hat f(\overline u)\] 
for all $t$ and $x\in\mathbb R$. 
\hfill 
$\blacksquare$

%
%
%
\newpage

\noindent 
{\bf Proof of Theorem 1.} 

When $\mu(\{0\})=1$, the conclusion of Theorem 1 is trivial. Suppose $\mu(\{0\})\not=1$. 
Then, in virtue of Lemmas 10, 11 and $f^\prime(\alpha)>0$, 
the semiflow $\{Q^t\}_{t\in[0,+\infty)}$ 
on $\mathcal M$ generated by (1.1) satisfies Hypothesis 7. So, Theorem 8 can work. 

We take a Lipschitz continuous function $\hat f$ on $\mathbb R$ with (4.19).  
Then, by Lemma 12, there exist two constants $\underline c$, $\overline c$, 
two bounded, continuous and monotone functions $\hat {\underline \psi}$ 
and $\hat {\overline \psi}$ 
on $\mathbb R$ with 
$\hat {\underline \psi}(0)\in(-\infty,0)$, $\hat {\underline \psi}(+\infty)\in(\alpha,1)$, 
$\hat {\overline \psi}(-\infty)\in(0,\alpha)$ and $\hat {\overline \psi}(0)\in(1,+\infty)$ 
such that ${\underline u}(t,x):=\hat {\underline \psi}(x+\underline ct)$ 
is a sub-solution to (4.20) 
and ${\overline u}(t,x):=\hat {\overline \psi}(x+\overline ct)$ 
is a super-solution to (4.20). 

We put $\underline \psi:=\max\{\hat {\underline \psi},0\}\in\mathcal M$ 
and $\overline \psi:=\min\{\hat {\overline \psi},1\}\in\mathcal M$. 
Then, ${\underline \psi}(0)=0$, 
${\underline \psi}(+\infty)\in(\alpha,1)$, 
${\overline \psi}(-\infty)\in(0,\alpha)$ 
and ${\overline \psi}(0)=1$ hold. Further, $(Q^t[\underline \psi])(x)$ 
and $(Q^t[\overline \psi])(x)$ are solutions to not only (1.1) 
but also (4.20) in $t\in[0,+\infty)$. 
Hence, because $\hat {\underline \psi}\leq\underline \psi$ 
and $\overline \psi\leq\hat {\overline \psi}$ hold, 
we have 
\[{\underline \psi}(x+\underline ct)
=\max\{\hat {\underline \psi}(x+\underline ct),0\}
\leq(Q^t[\underline \psi])(x)\]
and 
\[(Q^t[\overline \psi])(x)
\leq\min\{\hat {\overline \psi}(x+\overline ct),1\}
={\overline \psi}(x+\overline ct)\] 
for all $t\in[0,+\infty)$. Therefore, by Theorem 8, 
there exist $c\in[\underline c,\overline c]$ 
and $\phi\in\mathcal M$ with $\phi(-\infty)=0$ 
and $\phi(+\infty)=1$ such that $(Q^t[\phi])(x-ct)\equiv \phi(x)$ holds 
for all $t\in[0,+\infty)$. So, $u(t,x):=\phi(x+ct)$ is a solution to (1.1). 
\hfill 
$\blacksquare$

\sect{Appendix} 
\( \, \, \, \, \, \, \, \) 
In this section, we recall some known results from [17]. 
We use them in Section 3 to prove the abstract theorems 
stated in Section 2. 


The following is the same as Proposition 10 of [17]: 
\begin{proposition} \ 
Let a map $Q_0:\mathcal M \rightarrow \mathcal M$ 
satisfy {\rm Hypotheses 2} {\rm (i)}, {\rm (ii)} and {\rm (iii)}. 
Suppose a sequence $\{u_k\}_{k\in \mathbb N} \subset \mathcal M$ 
converges to $u\in \mathcal M$ almost everywhere. 
Then, $\lim_{k\rightarrow \infty}(Q_0[u_k])(x)=(Q_0[u])(x)$ holds 
for all continuous points $x\in\mathbb R$ of $Q_0[u]$. 
\end{proposition}


The following is the same as Lemma 11 of [17]: 
\begin{lemma} \ 
Let a sequence $\{u_k\}_{k\in \mathbb N}$ of monotone nondecreasing functions on $\mathbb R$ 
converge to a monotone nondecreasing function $u$ on $\mathbb R$ almost everywhere. 
Then, $\lim_{k\rightarrow \infty}u_k(x_k)=u(x)$ holds 
for all continuous points $x\in\mathbb R$ of $u$
and sequences $\{x_k\}_{k\in \mathbb N} \subset \mathbb R$ 
with $\lim_{k\rightarrow \infty}x_k=x$. 
\end{lemma}


The following is the same as Lemma 12 of [17]: 
\begin{lemma} \ 
Let $Q=\{Q^t\}_{t\in[0,+\infty)}$ be a family of maps from $\mathcal M$ to $\mathcal M$. 
Suppose $Q$ satisfies {\rm Hypothesis 6 (ii)}. 
Then, $\lim_{t\rightarrow 0}(Q^t[u])(x-ct)=u(x)$ holds 
for all $u\in \mathcal M$, $c\in\mathbb R$ 
and continuous points $x\in\mathbb R$ of $u$. 
\end{lemma}


The following is the same as Theorem 6 of [17] (Theorem 5 of [18]): 
\begin{lemma} \ 
Let $Q^t$ be a map from $\mathcal M$ to $\mathcal M$ for $t\in[0,+\infty)$. 
Suppose the map $Q^t$ satisfies  {\rm Hypotheses 2 (i)}, {\rm (ii)} and {\rm (iii)} 
for all $t\in (0,+\infty)$, and the family $Q:=\{Q^t\}_{t\in[0,+\infty)}$ {\rm Hypotheses 6}. 
Suppose the map $Q^t$ is monostable for all $t\in (0,+\infty)$; i.e., 
\[0<\gamma<1 \Longrightarrow \gamma < Q^t[\gamma]\] 
for all $t\in (0,+\infty)$ and constant functions $\gamma$. 
Then, the following holds {\rm :} 

Let $c\in\mathbb R$. Suppose there exist $\tau\in(0,+\infty)$ and $\phi\in\mathcal M$ 
with $(Q^\tau[\phi])(x-c\tau)\leq \phi(x)$, $\phi\not\equiv 0$ and $\phi\not\equiv 1$. 
Then, there exists $\varphi\in\mathcal M$ with $\varphi(-\infty)=0$ and $\varphi(+\infty)=1$ 
such that $(Q^t[\varphi])(x-ct)\equiv \varphi(x)$ holds for all $t\in[0,+\infty)$. 
\end{lemma}

The following follows from Lemma 16: 
\begin{lemma} \ 
Let $Q^t$ be a map from $\mathcal M$ to $\mathcal M$ for $t\in[0,+\infty)$. 
Suppose the map $Q^t$ satisfies {\rm Hypotheses 2} for all $t\in (0,+\infty)$, 
and the family $Q:=\{Q^t\}_{t\in[0,+\infty)}$ {\rm Hypotheses 6}. 
Then, the following two hold {\rm :} 

{\rm (i)} \ Let $\tau\in(0,+\infty)$ and $c_-\in\mathbb R$. Suppose 
there exists $\phi_-\in\mathcal M$ with $(Q^\tau[\phi_-])(x-c_-\tau)\equiv \phi_-(x)$, 
$\phi_-(-\infty)=0$ and $\phi_-(+\infty)=\alpha$. Then, 
there exists $\varphi_-\in\mathcal M$ with $\varphi_-(-\infty)=0$ 
and $\varphi_-(+\infty)=\alpha$ such that 
$(Q^t[\varphi_-])(x-c_-t)\equiv \varphi_-(x)$ holds for all $t\in[0,+\infty)$. 

{\rm (ii)} \ Let $\tau\in(0,+\infty)$ and $c_+\in\mathbb R$. Suppose 
there exists $\phi_+\in\mathcal M$ with $(Q^\tau[\phi_+])(x-c_+\tau)\equiv \phi_+(x)$, 
$\phi_+(-\infty)=\alpha$ and $\phi_+(+\infty)=1$. Then, 
there exists $\varphi_+\in\mathcal M$ with $\varphi_+(-\infty)=\alpha$ 
and $\varphi_+(+\infty)=1$ such that 
$(Q^t[\varphi_+])(x-c_+t)\equiv \varphi_+(x)$ holds for all $t\in[0,+\infty)$. 
\end{lemma}
{\bf Proof.} We show (i). 
Put a set as 
\[\mathcal M_-:=\{u\, |\, u \text{ is a monotone nondecreasing}\]
\[\text{ and left continuous function on } \mathbb R \text{ with } 0\leq u\leq \alpha\}.\] 
Put two maps $R_- : \mathcal M \rightarrow \mathcal M_-$ 
and $S_- : \mathcal M_- \rightarrow \mathcal M$ as 
\[(R_-[u])(x):=\alpha(1-\lim_{h\downarrow+0}u(-x+h))\]
and 
\[(S_-[u_-])(x):=-\frac{1}{\alpha}(\lim_{h\downarrow+0}u_-(-x+h))+1.\]
So, the maps $R_-$ and $S_-$ are inverse in each other. 
In virtue of Hypotheses 2 (ii) and (iv), we can define 
a map $Q_-^t:\mathcal M \rightarrow \mathcal M$ by 
\[Q_-^t:=S_-\circ Q^t\circ R_-\]
for $t\in[0,+\infty)$. Then, $Q_-:=\{Q_-^t\}_{t\in[0,+\infty)}$ 
satisfies the assumption of Lemma 16. 
Hence, Lemma 16 works for the semiflow $Q_-$. 
Let ${\tilde \phi}_-:=S_-[\phi_-]\in\mathcal M$. Then, 
$(Q_-^\tau[{\tilde \phi}_-])(x-c_-\tau)\equiv {\tilde \phi}_-(x)$, 
${\tilde \phi}_-(-\infty)=0$ and ${\tilde \phi}_-(+\infty)=1$ hold. 
Therefore, by Lemma 16, there exists ${\tilde \varphi}_-\in\mathcal M$ 
with ${\tilde \varphi}_-(-\infty)=0$ and ${\tilde \varphi}_-(+\infty)=1$ 
such that $(Q_-^t[{\tilde \varphi}_-])(x-c_-t)\equiv {\tilde \varphi}_-(x)$ 
holds for all $t\in[0,+\infty)$. 
So, as we put $\varphi_-:=R_-[{\tilde \varphi}_-]\in\mathcal M_-$, 
we obtain the conclusion of (i). 

We show (ii). Put a set as 
\[\mathcal M_+:=\{u\, |\, u \text{ is a monotone nondecreasing}\]
\[\text{ and left continuous function on } \mathbb R \text{ with } \alpha\leq u\leq 1\}.\] 
Put two maps $R_+ : \mathcal M \rightarrow \mathcal M_+$ 
and $S_+ : \mathcal M_+ \rightarrow \mathcal M$ as 
\[(R_+[u])(x):=(1-\alpha)u(x)+\alpha\]
and 
\[(S_+[u_+])(x):=\frac{1}{1-\alpha}(u_+(x)-\alpha).\] 
Then, almost similarly as (i), we can obtain the conclusion of (ii). 
\hfill 
$\blacksquare$

\vspace*{1.6em} 

\noindent Acknowledgments. \ 
I am grateful to the anonymous
referees for comments. 
It was partially supported 
by Grant-in-Aid for Scientific Research (No.19740092) 
from Ministry of Education, Culture, Sports, Science 
and Technology, Japan. 

\newpage

\[ \begin{array}{c} \mbox{R\scriptsize EFERENCES}  \end{array} \]

[1] N. D. Alikakos, P. W. Bates and X. Chen, 
Periodic traveling waves and locating oscillating patterns in multidimensional domains, 
{\it Trans. Amer. Math. Soc.}, 351 (1999), 2777-2805. 

[2] P. Bates and F. Chen, 
Periodic traveling waves for a nonlocal integro-differential model, 
{\it Electron. J. Diff. Eqns.}, 1999 (1999), No. 26. 

[3] P. W. Bates and A. Chmaj, A discrete convolution model for phase transitions,  
{\it Arch. Rational Mech. Anal.}, 150 (1999), 281-305. 

[4] P. W. Bates, P. C. Fife, X. Ren and X. Wang, Traveling waves in a convolution model 
for phase transitions, {\it Arch. Rational Mech. Anal.}, 138 (1997), 105-136. 

[5] X. Chen, Existence, uniqueness, and asymptotic stability of traveling waves 
in nonlocal evolution equations, {\it Adv. Differential Equations}, 2 (1997), 125-160. 

[6] X. Chen, J.-S. Guo and C.-C. Wu, 
Traveling waves in discrete periodic media for bistable dynamics, 
{\it Arch. Rational Mech. Anal.}, 189 (2008), 189-236. 

[7] A. J. J. Chmaj, Existence of traveling waves for the nonlocal Burgers equation, 
{\it Appl. Math. Lett.}, 20 (2007), 439-444.  

[8] S.-N. Chow, J. Mallet-Paret and W. Shen, 
Traveling waves in lattice dynamical systems, 
{\it J. Differential Equations}, 
149 (1998), 248-291.  

[9] R. Coutinho and B. Fernandez, Fronts in extended systems of bistable maps 
coupled via convolutions, {\it Nonlinearity}, 17 (2004), 23-47. 

[10] J. Coville, On uniqueness and monotonicity of solutions 
of non-local reaction diffusion equation, {\it Ann. Mat. Pura Appl.},  
185 (2006), 461-485. 

[11] J. Coville, Travelling fronts in asymmetric nonlocal reaction diffusion equation: 
The bistable and ignition case, preprint. 

[12] P. C. Fife and J. B. McLeod, A phase plane discussion of convergence to travelling fronts 
for nonlinear diffusion, {\it Arch. Rational Mech. Anal.}, 75 (1981), 281-314. 

[13] F. Hamel and S. Omrani, Existence of multidimensional travelling fronts 
with a multistable nonlinearity, 
{\it Adv. Differential Equations}, 5 (2000), 557-582. 

[14] J. Mallet-Paret, The global structure of traveling waves 
in spatially discrete dynamical systems, {\it J. Dynam. Differential Equations},  
11 (1999), 49-127. 

[15] A. Volpert and V. Volpert, Existence of multidimensional travelling waves and systems of waves, 
{\it Comm. Partial Differential Equations}, 26 (2001), 421--459.

[16] H. F. Weinberger, Long-time behavior of a class of biological models, 
{\it SIAM J. Math. Anal.}, 13 (1982), 353-396. 

[17] H. Yagisita, 
Existence of traveling wave solutions 
for a nonlocal monostable equation: an abstract approach, 
{\it Discrete Contin. Dyn. Syst.}, submitted 
({\bf For editors and referees:} http://arxiv.org/abs/0807.3612). 

[18] H. Yagisita, Existence and nonexistence of traveling waves 
for a nonlocal monostable equation, {\it Publ. Res. Inst. Math. Sci.}, submitted.


\end{document}